\documentclass[a4paper,10pt]{article}

\usepackage[utf8]{inputenc}
\usepackage{amsmath,amssymb}
\newtheorem{thm}{Theorem}[section]
\newtheorem{lem}[thm]{Lemma}
\newtheorem{pro}[thm]{Proposition}
\newtheorem{cor}[thm]{Corollary}
\usepackage[all]{xy}

\title{\bf On Gieseker stability for Higgs sheaves}
\author{S. A. H. Cardona and O. Mata-Guti\'errez}

\begin{document}

\maketitle

\begin{abstract}
We review the notion of Gieseker stability for torsion-free Higgs sheaves. This notion is a natural generalization of the classical notion of Gieseker stability for torsion-free coherent sheaves. In this article we prove some basic properties that are similar to the classical ones for torsion-free coherent sheaves over projective algebraic manifolds. In particular, we show that Gieseker stability for torsion-free Higgs sheaves can be defined using only Higgs subsheaves with torsion-free quotients; we also prove that a direct sum of two Higgs sheaves is Gieseker semistable if and only if the Higgs sheaves are both Gieseker semistable with equal normalized Hilbert polynomial; then we prove that a classical property of morphisms between Gieseker semistable sheaves also holds in the Higgs case; as a consequence of this and because of an existing relation between Mumford-Takemoto stability and Gieseker stability for Higgs sheaves, we obtain certain properties concerning the existence of Hermitian-Yang-Mills metrics, simplesness and extensions. Finally, we make some comments about Jordan-H\"older and Harder-Narasimhan filtrations for Higgs sheaves. \\


\noindent{\it MS Classification}: 53C07, 53C55, 32C15.

\end{abstract}

\section{Introduction}

As it is well known, Mumford \cite{Mumford} introduced a notion of stability for vector bundles over curves and such a notion was latter on generalized to algebraic surfaces by Takemoto \cite{Takemoto 1, Takemoto 2}, who called this 
$H$-stability, where $H$ denoted an ample line bundle over the base manifold. This notion was also studied by Kobayashi \cite{Kobayashi} when the base manifold was a compact K\"ahler manifold, he called it Mumford-Takemoto stability. In this article we will refer to Mumford-Takemoto stability simply as stability. \\

On the other hand, the notion of Gieseker $H$-stability was introduced by Gieseker \cite{Gieseker} and 
Kobayashi \cite{Kobayashi}. In their works, this kind of stability was studied for torsion-free coherent sheaves over projective algebraic manifolds. Here again, $H$ denoted an ample line bundle over the base manifold. In this article we will refer to Gieseker $H$-stability simply as Gieseker 
stability. Now, Kobayashi proved some results concerning Gieseker stability and its relation with stability. 
In particular, he showed that stability implies Gieseker stability, and Gieseker semistability implies 
semistability. He proved also a relation between Gieseker stability and the existence of 
Hermitian-Yang-Mills metrics (from now on abbreviated as HYM-metrics) on holomorphic vector bundles over projective algebraic manifolds. Indeed, as we will see, by 
using a famous theorem of Bando and Siu \cite{Bando-Siu}, this result can be extended to reflexive sheaves, where 
the differential geometric counterpart in this case is the notion of admissible 
HYM-metric.\footnote{A hermitian metric on a coherent sheaf is one defined on the locally 
free part of it; such a metric is said to be {\it admissible}, if the Chern curvature is square integrable and 
its trace with respect to the K\"ahler form is uniformly bounded (see \cite{Bando-Siu, Biswas-Schumacher} for 
more details).} Additionally, Kobayashi proved the injectivity or surjectivity of any morphism between Gieseker semistable 
sheaves, depending on relations between certain invariants of the sheaves (their normalized Hilbert 
polynomials) and the Gieseker stability of one of the sheaves.\\

Now, Higgs bundles were first introduced by Hitchin \cite{Hitchin} as geometric objects associated 
to solutions of the selfdual Yang-Mills equations over curves. Later on, Simpson \cite{Simpson} generalized the 
ideas of Hitchin to base manifolds of higher dimensions. In his article, Simpson defines Higgs bundles over 
compact (and some non-compact) K\"ahler manifolds, and proves a Hitchin-Kobayashi correspondence for such objects. 
Namely, he proves that a Higgs bundle admits a HYM-metric if and only if it is (Mumford-Takemoto) polystable. Since such a 
correspondence makes reference to the notion of stability, it was necessary to consider in the theory 
Higgs sheaves and not only Higgs bundles. Now, Simpson \cite{Simpson M1, Simpson M2} introduced the notion of 
Gieseker stability for a kind of more general objects, known as $\Lambda$-modules. Higgs sheaves are particular cases of 
$\Lambda$-modules, and hence some (but not all) properties of Higgs sheaves can be obtained directly from properties of 
$\Lambda$-modules. Indeed, since the category of Higgs sheaves is also an abelian category (see \cite{Cardona 2} for details), 
it is possible to obtain other properties for these objects using this fact. \\
  
On the other hand, several properties of holomorphic bundles (resp. coherent sheaves) have been extended to Higgs bundles (resp. Higgs sheaves). In particular, 
there are Bochner's vanishing theorems for Higgs bundles \cite{Cardona 3}; and there are extensions of the Hitchin-Kobayashi 
correspondence to polystable reflexive Higgs sheaves \cite{Biswas-Schumacher} (a reflexive Higgs sheaf is polystable if and 
only if it has an admissible HYM-metric) and to semistable Higgs bundles \cite{Bruzzo-Granha, Jiayu-Zhang}; in 
this case, the differential geometric counterpart is the notion of 
approximate Hermitian-Yang-Mills metric (from now on abbreviated apHYM-metric). More recently \cite{Cardona 4}, 
the notion of $T$-stability (also called Bogomolov stability) has been studied in the context of Higgs sheaves over 
compact complex manifolds. In this article, we study the notion of Gieseker stability for Higgs sheaves and we prove some basic properties, as we said before, some of these results can be seen as natural extensions of classical results for coherent sheaves.\\

This article is organized as follows, in the second section we review some basics notions concerning torsion-free 
Higgs sheaves and the theory of Chern characters. In particular, we define the normalized Hilbert polynomial and we observe that a classical identity involving ranks and these polynomials, also 
applies to Higgs sheaves. In the 
third section, we introduce the notion of Gieseker stability for torsion-free Higgs sheaves as a natural generalization 
of the classical notion of Gieseker stability for torsion-free coherent sheaves over projective algebraic 
manifolds and then we prove similar results to the classical ones for Higgs sheaves. In particular, we show that the existing 
relation between Gieseker stability and stability extends naturally to Higgs sheaves. Then, we show 
that in analogy to stability, the Gieseker stability can be defined using only Higgs subsheaves 
with torsion-free quotients. Next, we prove that a direct sum of two Higgs sheaves is Gieseker semistable if 
and only if the Higgs sheaves are both Gieseker semistable with equal normalized Hilbert polynomial. At the 
end of this section we prove our main result. Namely, we prove the vanishing, injectivity or surjectivity of any morphism between Gieseker 
semistable Higgs sheaves, depending on certain relations between their normalized Hilbert polynomials and 
(eventually) the Gieseker stability of one of these sheaves. Finally, by using the Hitchin-Kobayashi correspondence for reflexive 
Higgs sheaves \cite{Biswas-Schumacher}, we prove a relation between the Gieseker stability and the 
existence of admissible Hermitian-Yang-Mills metrics for reflexive Higgs sheaves. In the final section and by following the standard literature \cite{Gieseker, Kobayashi}, we review the notions of Jordan-H\"older and Harder-Narasimhan filtrations for Higgs sheaves. It is important to mention that the existence of filtrations of these type has been already studied by Simpson \cite{Simpson M1, Simpson M2} in a more general setting. Indeed, Simpson proves the existence of Jordan-H\"older and Harder-Narasimhan filtrations for 
$\Lambda$-modules, since Higgs sheaves are particular cases of these objects, we really know that there exist Jordan-H\"older and Harder-Narasimhan filtrations for Higgs sheaves.

\section{Preliminaries}
In this section we fix the notation and review some basic definitions in the projective case. Let $X$ be an 
$n$-dimensional projective algebraic manifold and let $\Omega_{X}^{1}$ be the cotangent sheaf to $X$, i.e., 
it is the sheaf of holomorphic one-forms on $X$. Then a Higgs sheaf ${\mathfrak E}$ over $X$ is a pair 
$(E,\phi)$, where $E$ is a coherent sheaf over $X$ and $\phi : E\rightarrow E\otimes\Omega_{X}^{1}$ is a morphism 
of ${\cal O}_{X}$-modules (usually called the Higgs field) such that $\phi\wedge\phi : E\rightarrow E\otimes\Omega_{X}^{2}$ 
vanishes. Following \cite{Bruzzo-Granha}, a section $s$ of $E$ such that $\phi(s)=s\otimes\lambda$ for some holomorphic one-form $\lambda$ on $X$, is called a $\phi$-invariant section of ${\mathfrak E}$. A Higgs sheaf ${\mathfrak E}$ is said to be torsion-free (resp. locally free, reflexive, torsion) if the coherent sheaf $E$ is torsion-free 
(resp. locally free, reflexive, torsion); and by definition a Higgs bundle is a locally free Higgs sheaf. Also, 
the singularity set of a Higgs sheaf ${\mathfrak E}$ is defined as the singularity set of the corresponding 
coherent sheaf $E$ and similarly for the rank and first Chern class, i.e., $S({\mathfrak E})=S(E)$, ${\rm rk\,}{\mathfrak E}={\rm rk\,}E$ and 
$c_{1}({\mathfrak E})=c_{1}(E)$. \\

A Higgs subsheaf ${\mathfrak F}$ of ${\mathfrak E}$ is a subsheaf $F$ of $E$ such that 
$\phi(F)\subset F\otimes\Omega_{X}^{1}$, and hence the pair ${\mathfrak F}=(F,\phi|_{F})$ becomes itself a Higgs 
sheaf. In this case, 
the quotient $E/F$ has a natural Higgs morphism induced from the Higgs morphisms of ${\mathfrak F}$ and 
${\mathfrak E}$ and hence it is a Higgs sheaf, that we denote by ${\mathfrak E}/{\mathfrak F}$. Now, if ${\mathfrak E}$ and ${\mathfrak E}'$ 
are two Higgs sheaves over $X$, a morphism $f:{\mathfrak E}\longrightarrow{\mathfrak E}'$ is a map 
$f:E\longrightarrow E'$ between the corresponding coherent sheaves such that 
the diagram 
\begin{displaymath}
 \xymatrix{
E \ar[r]^{\phi} \ar[d]^{f}    &     E\otimes\Omega_{X}^{1} \ar[d]^{f\otimes 1}   \\
E' \ar[r]^{\phi'}              &     E'\otimes\Omega_{X}^{1}   \\
}
\end{displaymath}
commutes, in particular an endomorphism $f$ of ${\mathfrak E}$ is just a morphism $f:{\mathfrak E}\longrightarrow{\mathfrak E}$. The kernel and image of any morphism of Higgs sheaves are again Higgs sheaves in a natural way, and the 
torsion subsheaf of any Higgs sheaf is again a Higgs sheaf. An exact sequence of Higgs sheaves is just an exact sequence of their corresponding coherent sheaves, 
such that each morphism is a morphism of Higgs sheaves; in particular, any Higgs subsheaf ${\mathfrak F}$ of 
${\mathfrak E}$ defines a short exact sequence of Higgs sheaves 
\begin{equation}
 \xymatrix{
0 \ar[r]  &  {\mathfrak F} \ar[r]  &   {\mathfrak E} \ar[r]  &   {\mathfrak Q} \ar[r] &  0  
} \label{Higgs extension}
\end{equation}
where ${\mathfrak Q}={\mathfrak E}/{\mathfrak F}$. A short exact sequence of this type is usually called a Higgs 
extension (for more details on all of these basic properties and definitions see for instance \cite{Biswas-Schumacher, Bruzzo-Granha} or 
\cite{Cardona 2}).\\

Let ${\mathfrak E}$ be a torsion-free Higgs sheaf, as it was shown in \cite{Cardona 4} we can use a canonical isomorphism of coherent sheaves to construct a Higgs morphism for ${\rm det\,}E$, and hence, this determinant can be considered as a Higgs line bundle ${\rm det\,}{\mathfrak E}$. Now, 
from a classical result of coherent sheaves \cite{Kobayashi}, we know that any injective morphism between torsion-free sheaves of the same rank induces an injective morphism between its determinant bundles. From all this we get the following result.
\begin{pro}\label{prop. morph. det}
Let ${\mathfrak E}$ and ${\mathfrak E}'$ be two torsion-free Higgs sheaves over $X$. If ${\rm rk\,}{\mathfrak E}={\rm rk\,}{\mathfrak E}'$ and $f:{\mathfrak E}\longrightarrow{\mathfrak E}'$ is an injective morphism, then the induced morphism ${\rm det\,}f:{\rm det\,}{\mathfrak E}\longrightarrow{\rm det\,}{\mathfrak E}'$ is also an injective morphism.
\end{pro}

On the other hand, since $X$ is a projective algebraic manifold, we can consider an ample line bundle over it. Let 
$H$ be a fixed ample line bundle over $X$, then we will define the notions of stability and Gieseker stability 
with respect to $H$. The ample line bundle $H$ can be considered as a Higgs line bundle with 
zero Higgs field\footnote{Even when ${\mathfrak H}$ is in essence $H$, we introduce this 
definition in order to simplify our notation, e.g., in that way we can do tensor products of any Higgs sheaf 
${\mathfrak E}$ with powers of ${\mathfrak H}$ and the result is considered as a Higgs sheaf.} ${\mathfrak H}=(H,0)$. Now, let $\omega_{H}$ be a representative of $c_{1}({\mathfrak H})$, then as it is well 
known \cite{Bruzzo-Granha, Simpson} the degree of a Higgs sheaf ${\mathfrak E}$ can be defined with respect 
to this $\omega_{H}$ as
\begin{equation}
 {\rm deg\,}{\mathfrak E} = \int_{X}c_{1}(\mathfrak E)\cdot\omega_{H}^{n-1}\, \nonumber
\end{equation}
where the dot here (and from now on) is the usual notation in literature for the wedge product of forms. The slope of ${\mathfrak E}$ is 
given by $\mu_{\mathfrak E}={\rm deg\,}{\mathfrak E}/{\rm rk\,}{\mathfrak E}$ and we say that a torsion-free Higgs sheaf ${\mathfrak E}$ 
over $X$ is stable (resp. semistable) if for any Higgs subsheaf ${\mathfrak F}$ of ${\mathfrak E}$ with 
$0<{\rm rk\,}{\mathfrak F}<{\rm rk\,}{\mathfrak E}$ we have $\mu_{\mathfrak F}<\mu_{\mathfrak E}$ (resp. $\le$).\\

As it is well known, Hirzebruch \cite{Hirzebruch} defines Chern classes and proves a Riemann-Roch formula for 
holomorphic bundles over projective algebraic manifolds. Now, the theory of Chern classes and characters has 
been extended to coherent sheaves by O'Brian, Toledo and Tong \cite{Toledo 1, Toledo 2} and they prove 
a Hirzebruch-Riemann-Roch formula in that case. Applied to Higgs sheaves, this formula says that for a Higgs sheaf 
${\mathfrak E}$ over $X$ 
\begin{equation}
 \sum_{i\geq 0}(-1)^{i}{\rm dim\,}H^{i}(X,{\mathfrak E}) = \int_{X}{\rm ch}({\mathfrak E})\cdot{\rm td}(X)\, \label{Hirzebruch-Riemann-Roch}
\end{equation}
where $H^{i}(X,{\mathfrak E})$ denotes the $i$-th cohomology group of $X$ with coefficients in 
${\mathfrak E}$, and ${\rm ch}({\mathfrak E})$ and ${\rm td}(X)$ denote the Chern character and the Todd class 
of $X$, respectively. The left hand side of (\ref{Hirzebruch-Riemann-Roch}) is commonly called the Euler-Poincar\'e 
characteristic of ${\mathfrak E}$ and is denoted by $\chi(X,{\mathfrak E})$ or simply by $\chi({\mathfrak E})$. 
Similar to the slope $\mu$, the Euler-Poincar\'e characteristic $\chi$ satisfies some elementary properties; for 
instance, since the Higgs extension (\ref{Higgs extension}) is in particular a short exact sequence of coherent 
sheaves, we have (see \cite{Hirzebruch} for a proof of this)
\begin{equation}
 \chi(\mathfrak E) = \chi(\mathfrak F) + \chi(\mathfrak Q)\,.  \label{Basic prop chi}
\end{equation}

If we consider again the Higgs line bundle ${\mathfrak H}$ that we have introduced before, then associated with any    
torsion-free Higgs sheaf ${\mathfrak E}$ of positive rank, we have a (torsion-free) Higgs sheaf 
${\mathfrak E}(k)={\mathfrak E}\otimes{\mathfrak H}^{k}$ for $k\in {\mathbb Z}$ and we define 
\begin{equation}
p_{\mathfrak E}(k)=\chi({\mathfrak E}(k))/{\rm rk\,}{\mathfrak E}\,.\label{Reduced Hilbert polynomial}\end{equation}
In the literature \cite{Hartshorne, Huybrechts-Lehn} the Euler-Poincar\'e characteristic of ${\mathfrak E}(k)$, $\chi({\mathfrak E}(k))$, is commonly called the Hilbert polynomial of ${\mathfrak E}$, and Simpson \cite{Simpson M1} called the quotient 
$p_{\mathfrak E}(k)$ the {\it normalized Hilbert polynomial} of ${\mathfrak E}$. Now, from basic properties of Chern characters we know 
that ${\rm ch}({\mathfrak E}\otimes{\mathfrak H}^{k}) = {\rm ch}({\mathfrak E})\,{\rm ch}({\mathfrak H}^{k})$, and hence 
the Hirzebruch-Riemann-Roch formula (\ref{Hirzebruch-Riemann-Roch}) for ${\mathfrak E}(k)$ becomes
\begin{equation}
 \chi({\mathfrak E}(k)) = \int_{X}{\rm ch}({\mathfrak E})\cdot{\rm ch}({\mathfrak H}^{k})\cdot{\rm td}(X)\,. \label{Hirzebruch-Riemann-Roch Ek}
\end{equation}
Now, from \cite{Kobayashi} we really have expressions for ${\rm ch}({\mathfrak E})$, ${\rm ch}({\mathfrak H}^{k})$ and 
${\rm td}(X)$ in terms of the Chern classes of ${\mathfrak E}$, ${\mathfrak H}$ and $X$, respectively. 
To be precise, if we set $r={\rm rk\,}{\mathfrak E}$ we have
\begin{equation}
 {\rm ch}({\mathfrak E}) = r + c_{1}({\mathfrak E}) + \frac{1}{2}(c_{1}({\mathfrak E})^{2} - 2c_{2}({\mathfrak E})) + \cdots\,, \nonumber
\end{equation}
\begin{equation}
 {\rm ch}({\mathfrak H}^{k}) = 1 + kc_{1}({\mathfrak H}) + \frac{1}{2}k^{2}c_{1}({\mathfrak H})^{2} + \cdots +  \frac{1}{n!}k^{n}c_{1}({\mathfrak H})^{n}\,, \nonumber
\end{equation}
\begin{equation}
 {\rm td}(X) = 1 +  \frac{1}{2}c_{1}(X) + \frac{1}{12}(c_{1}(X)^{2} + c_{2}(X)) + \cdots\,. \nonumber
\end{equation}
If we replace these expressions for Chern characters in (\ref{Hirzebruch-Riemann-Roch Ek}) we obtain 
\begin{equation}
 \chi({\mathfrak E}(k)) = \int_{X}\frac{r k^{n}}{n!}c_{1}({\mathfrak H})^{n} + \int_{X}\frac{k^{n-1}}{(n-1)!}c_{1}({\mathfrak H})^{n-1}
                           \cdot(c_{1}(\mathfrak E) + \frac{r}{2}c_{1}(X)) + \cdots \,, \label{chi Ek}
\end{equation}
where the dots represent all terms of lower order than $k^{n-1}$. Now, by dividing the expression (\ref{chi Ek}) 
by $r$ and using $\omega_{H}$ to define a volume form of $X$ and a degree of ${\mathfrak E}$, we get the 
formula
\begin{equation}
 p_{\mathfrak E}(k) = k^{n}{\rm vol\,}X + \frac{k^{n-1}}{(n-1)!}\mu_{\mathfrak E}  
                        + \frac{k^{n-1}}{2(n-1)!}\int_{X}c_{1}(X)\cdot{\omega_{H}}^{n-1} + \cdots \,, 
 \label{p Ek}
\end{equation}
where the dots again represent terms of lower order than $k^{n-1}$. From formula (\ref{p Ek}) it follows that, 
up to multiplication by the factor ${\rm deg\,}{\mathfrak H}$ (which depends on $\omega_{H}$), the reduced Hilbert 
polynomial in \cite{Huybrechts-Lehn} is in essence the same $p_{\mathfrak E}(k)$. Now, notice that in analogy to the slope $\mu$,
the polynomial $p_{\mathfrak E}(k)$ for Higgs sheaves is the same as the classical polynomial $p(E(k))$ for coherent sheaves of Kobayashi in \cite{Kobayashi}; and as it is well known, there exists a result involving these classical polynomials and ranks for any short exact sequence; hence such a result applies directly to Higgs sheaves and can be written as follows.    
\begin{lem}\label{Rel. p-rk}
 Let us consider the Higgs extension (\ref{Higgs extension}) over $X$. Then for any integer $k$ we have 
 \begin{equation}
  {\rm rk\,}{\mathfrak F}\,(p_{\mathfrak E}(k)-p_{\mathfrak F}(k)) + {\rm rk\,}{\mathfrak Q}\,(p_{\mathfrak E}(k)-p_{\mathfrak Q}(k)) = 0\,. \nonumber
 \end{equation}
\end{lem}
\noindent {\it Proof:} Let $r,r'$ and $r''$ be the ranks of ${\mathfrak E},{\mathfrak F}$ and ${\mathfrak Q}$, respectively. 
Since ${\mathfrak H}$ is a Higgs line bundle, tensoring (\ref{Higgs extension}) by ${\mathfrak H}^{k}$ we have the Higgs extension
\begin{equation}
 \xymatrix{
0 \ar[r]  &  {\mathfrak F}(k) \ar[r]  &   {\mathfrak E}(k) \ar[r]  &   {\mathfrak Q}(k) \ar[r] &  0  
} \nonumber
\end{equation}
and because $r=r' +r''$, by applying (\ref{Basic prop chi}) to this sequence we have
\begin{equation}
 (r' + r'')p_{\mathfrak E}(k) = \chi({\mathfrak E}(k)) = \chi({\mathfrak F}(k)) + \chi({\mathfrak Q}(k)) = 
 r'p_{\mathfrak F}(k) + r''p_{\mathfrak Q}(k) \nonumber
\end{equation}
and the result follows.  \;\;Q.E.D.  \\

Finally, we introduce the following notation of Gieseker \cite{Gieseker}. Let ${\mathfrak E}$ and ${\mathfrak E}'$ 
be two Higgs sheaves over $X$ and let $p_{\mathfrak E}$ and $p_{{\mathfrak E}'}$ denote the corresponding 
normalized Hilbert polynomials. Then, we say that $p_{\mathfrak E}\prec p_{{\mathfrak E}'}$ 
(resp. $\preceq$) if the inequality $p_{\mathfrak E}(k)<p_{{\mathfrak E}'}(k)$ (resp. $\le$) holds for 
sufficiently large integers $k$. Notice that from this definition $p_{\mathfrak E}\preceq p_{{\mathfrak E}'}$ if 
and only if $p_{\mathfrak E}\prec p_{{\mathfrak E}'}$ or $p_{\mathfrak E} = p_{{\mathfrak E}'}$, where the last 
expression means that the normalized Hilbert polynomials are equal as polynomials.\footnote{If two polynomials are equal for sufficiently large integers $k$, they are necessarily the same.} Notice also that $\prec$ is transitive, i.e., if ${\mathfrak E}, {\mathfrak E}'$ and ${\mathfrak E}''$ are Higgs sheaves over $X$ with $p_{\mathfrak E}\prec p_{{\mathfrak E}'}$ and $p_{{\mathfrak E}'}\prec p_{{\mathfrak E}''}$, then $p_{\mathfrak E}\prec p_{{\mathfrak E}''}$ (clearly a similar result holds for $\preceq$).

\section{Gieseker stability}
Let ${\mathfrak H}$ be the Higgs bundle defined in Section 2, we say that a torsion-free Higgs sheaf 
${\mathfrak E}$ over $X$ is Gieseker stable (resp. Gieseker semistable), if for every Higgs subsheaf 
${\mathfrak F}$ of ${\mathfrak E}$ with $0<{\rm rk\,}{\mathfrak F}<{\rm rk\,}{\mathfrak E}$, 
we have $p_{\mathfrak F}\prec p_{\mathfrak E}$ (resp. $\preceq$). We say that a Higgs sheaf is strictly Gieseker semistable, if it is Gieseker semistable but not Gieseker stable. From these definitions it is clear that any Higgs sheaf of rank one is Gieseker stable. Now, as a consequence 
of (\ref{p Ek}) we get a relation between the notions of stability and Gieseker stability, which is 
indeed an extension of a classical proposition in \cite{Kobayashi}. To be precise we have the following result.
\begin{pro}\label{H-s-->G-H-s}
 Let ${\mathfrak E}$ be a torsion-free Higgs sheaf over $X$. Then, \\
 {\bf (i)} If ${\mathfrak E}$ is stable, then it is Gieseker stable; \\
 {\bf (ii)} If ${\mathfrak E}$ is Gieseker semistable, then it is semistable.
\end{pro}
\noindent {\it Proof:} Assume first that ${\mathfrak E}$ is stable and let ${\mathfrak F}$ be a Higgs subsheaf of ${\mathfrak E}$ with 
 $0<{\rm rk\,}{\mathfrak F}<{\rm rk\,}{\mathfrak E}$. Then from (\ref{p Ek}) it follows that (for any integer $k$)
\begin{equation}
 p_{\mathfrak E}(k) - p_{\mathfrak F}(k) = \frac{k^{n-1}}{(n-1)!}(\mu_{\mathfrak E} - \mu_{\mathfrak F}) 
                                               + \cdots\,.  \label{p Ek - p Fk}
\end{equation}
Since $\mu_{\mathfrak F}<\mu_{\mathfrak E}$, for sufficiently large integers $k$ the term of order $k^{n-1}$ on 
the right hand side of (\ref{p Ek - p Fk}) becomes a dominant one; hence, the left hand side of (\ref{p Ek - p Fk}) becomes 
positive for such integers $k$ and $p_{\mathfrak F}\prec p_{\mathfrak E}$, which proves (i). Now, if ${\mathfrak E}$ is Gieseker semistable, 
then the left hand side of (\ref{p Ek - p Fk}) is non-negative for sufficiently large integers $k$, but this 
implies necessarily that $\mu_{\mathfrak F}\le\mu_{\mathfrak E}$, and hence (ii) follows.  \;\;Q.E.D.  \\  

On the other hand, from the polynomial expressions of Chern characters, it is easy to see that all terms of lower 
order than $k-1$ on the right hand side of formula (\ref{p Ek - p Fk}) contain either, higher dimensional Chern classes or 
products of first Chern classes. This shows that in the one-dimensional case these additional terms in 
(\ref{p Ek - p Fk}) are all zero and the notions of stability and Gieseker stability coincide.\footnote{Indeed, the 
same kind of argument shows that these notions coincide classically, i.e., for coherent sheaves over compact Riemann 
surfaces.} Now, since in the one-dimensional case there exist examples of stable Higgs bundles that are not stable in the classical sense 
(see \cite{Hitchin} for details), we know that the notion of Gieseker stability in the Higgs case is not the same 
classical Gieseker stability.\\

As it is well known \cite{Kobayashi}, the stability can be equivalently defined from an 
inequality involving quotients, instead of subsheaves. An analog result also holds for Gieseker stability 
and it extends straightforwardly to the Higgs case because of Lemma \ref{Rel. p-rk}. To be precise we have the following result.
\begin{pro}\label{Prop. G-Hss}
 Let ${\mathfrak E}$ be a torsion-free Higgs sheaf over $X$. Then, ${\mathfrak E}$ is 
 Gieseker stable (resp. Gieseker semistable) if and only if for every quotient Higgs sheaf ${\mathfrak Q}$ of 
 ${\mathfrak E}$ with $0<{\rm rk\,}{\mathfrak Q}<{\rm rk\,}{\mathfrak E}$, we have 
 $p_{\mathfrak E}\prec p_{\mathfrak Q}$ (resp. $\preceq$). 
\end{pro}

Now, as in the case of stability, the Gieseker stability can be defined without making reference 
to all Higgs subsheaves or Higgs quotients. This is in part a consequence of the formula (\ref{Basic prop chi}) 
and the fact that for $k$ sufficiently large, the term of order $k^{n}$ in (\ref{chi Ek}) becomes a dominant 
one. Hence we have the following result (its proof is similar to the proof in \cite{Cardona 2} of 
an analog result on stability).
\begin{pro}\label{Basic prop. torsion}
 Let ${\mathfrak E}$ be a torsion-free Higgs sheaf over $X$. Then, \\
 {\bf (i)} ${\mathfrak E}$ is Gieseker stable (resp. Gieseker semistable) if and only if for any 
 Higgs subsheaf ${\mathfrak F}$ of ${\mathfrak E}$ with $0<{\rm rk\,}{\mathfrak F}<{\rm rk\,}{\mathfrak E}$ and 
 ${\mathfrak E}/{\mathfrak F}$ torsion-free, we have $p_{\mathfrak F}\prec p_{\mathfrak E}$ 
 (resp. $\preceq$). \\
 {\bf (ii)} ${\mathfrak E}$ is Gieseker stable (resp. Gieseker semistable) if and only if for any 
 torsion-free Higgs quotients ${\mathfrak Q}$ of ${\mathfrak E}$ with $0<{\rm rk\,}{\mathfrak Q}<{\rm rk\,}{\mathfrak E}$, 
 we have $p_{\mathfrak E}\prec p_{\mathfrak Q}$ (resp. $\preceq$).
\end{pro}
\noindent {\it Proof:} There is nothing to prove in one direction. In order to prove the other one, suppose that the inequalities in (i) (resp. (ii)) hold for all proper Higgs subsheaves of ${\mathfrak E}$ of positive rank with torsion-free quotients (resp. for all torsion-free Higgs quotients of ${\mathfrak E}$ with rank strictly less than ${\rm rk\,}{\mathfrak E}$). \\

Let us consider the Higgs extension (\ref{Higgs extension}) of ${\mathfrak E}$, let ${\mathfrak T}$ be the 
torsion of ${\mathfrak Q}$ and ${\mathfrak Q}'={\mathfrak Q}/{\mathfrak T}$. Let ${\mathfrak F}'$ be the kernel 
of the natural composition ${\mathfrak E}\longrightarrow{\mathfrak Q}\longrightarrow{\mathfrak Q}'$. Then we obtain the following commutative diagram
\begin{displaymath}
 \xymatrix{
         &                                &                            &            0   \ar[d]           &   \\
         &           0  \ar[d]            &                            &  {\mathfrak T} \ar[d]           &   \\
0 \ar[r] & {\mathfrak F} \ar[r]\ar[d] & {\mathfrak E} \ar[r]\ar[d]^{\rm Id} & {\mathfrak Q} \ar[r]\ar[d] & 0 \\
0 \ar[r] & {\mathfrak F}'\ar[r]\ar[d] & {\mathfrak E} \ar[r]  & {\mathfrak Q}' \ar[r]\ar[d]        & 0 \\
         & {\mathfrak F}'/{\mathfrak F}\ar[d]  &                 &           0                     &   \\
         &                0               &                            &                                 &   \\
}
\end{displaymath}
where ${\mathfrak Q}'$ is torsion-free, ${\mathfrak F}$ is a Higgs subsheaf of ${\mathfrak F}'$ and 
${\mathfrak F}'/{\mathfrak F}\cong{\mathfrak T}$. Now, the vertical exact sequences in the above diagram are Higgs 
extensions for ${\mathfrak F}'$ and ${\mathfrak Q}$, tensoring these sequences by ${\mathfrak H}^{k}$ we get the 
following exact sequences 
\begin{equation}
 \xymatrix{
0 \ar[r]  &  {\mathfrak F}(k) \ar[r]  &   {\mathfrak F}'(k) \ar[r]  &   {\mathfrak T}(k) \ar[r] &  0  
} \label{Higgs extension A}
\end{equation}
\begin{equation}
 \xymatrix{
0 \ar[r]  &  {\mathfrak T}(k) \ar[r]  &   {\mathfrak Q}(k) \ar[r]  &   {\mathfrak Q}'(k) \ar[r] &  0  
} \label{Higgs extension B}
\end{equation}
and since for $k$ sufficiently large, the term of order $k^{n}$ in (\ref{chi Ek}) becomes a dominant one and its coefficient is positive (it is just $r{\rm Vol}\,X$), we can choose $\chi({\mathfrak T}(k))>0$. Then, by applying (\ref{Basic prop chi}) to the Higgs extensions (\ref{Higgs extension A}) and (\ref{Higgs extension B}) we obtain (for sufficiently large integers $k$) the inequalities
\begin{equation}
 \chi({\mathfrak F}'(k)) = \chi({\mathfrak F}(k)) + \chi({\mathfrak T}(k)) >  \chi({\mathfrak F}(k))\,, \nonumber
\end{equation}
\begin{equation}
 \chi({\mathfrak Q}(k)) = \chi({\mathfrak T}(k)) + \chi({\mathfrak Q}'(k)) >  \chi({\mathfrak Q}'(k))\,. \nonumber
\end{equation}
Since ${\mathfrak T}$ is torsion, ${\rm rk\,}{\mathfrak F}'={\rm rk\,}{\mathfrak F}$ and 
${\rm rk\,}{\mathfrak Q}={\rm rk\,}{\mathfrak Q}'$ and hence from the above inequalities we get  
\begin{equation} 
 p_{\mathfrak F} \prec p_{{\mathfrak F}'}\,, \quad\quad\quad p_{{\mathfrak Q}'} \prec p_{\mathfrak Q}\,. 
\label{inequalities}  
\end{equation}
On the other hand, since ${\mathfrak Q}'$  is torsion-free, from hypothesis we know that
\begin{equation} 
 p_{{\mathfrak F}'}\prec p_{\mathfrak E}\,, \quad\quad\quad p_{\mathfrak E}\prec p_{{\mathfrak Q}'}\,. 
 \label{inequalities 2}  
\end{equation}
At this point, from (\ref{inequalities}) and (\ref{inequalities 2}) it follows that ${\mathfrak E}$ is Gieseker stable. Now, if in the hypothesis we consider inequalities with $\preceq$, then in (\ref{inequalities 2}) we have $\preceq$ instead of $\prec$ and we obtain that ${\mathfrak E}$ is Gieseker semistable. \;\;Q.E.D.  \\

Notice that a torsion-free Higgs sheaf ${\mathfrak E}$ is not Gieseker stable if and only if there exists a Higgs subsheaf ${\mathfrak F}$ of it such that $p_{\mathfrak E}(k)\leq p_{\mathfrak F}(k)$ holds for sufficiently large integers $k$, i.e., if and only if $p_{\mathfrak E}\preceq p_{\mathfrak F}$ for some 
${\mathfrak F}$. Now, ${\mathfrak E}$ is Gieseker semistable if and only if for every Higgs subsheaf 
${\mathfrak F}$ the inequality $p_{\mathfrak F}(k)\leq p_{\mathfrak E}(k)$ holds for sufficiently large 
integers $k$. Therefore, ${\mathfrak E}$ is strictly Gieseker semistable if and only if there 
exists a Higgs subsheaf ${\mathfrak F}$ with $0<{\rm rk\,}{\mathfrak F}<{\rm rk\,}{\mathfrak E}$ 
such that $p_{\mathfrak F}=p_{\mathfrak E}$. From this fact and the proof of Proposition \ref{Basic prop. torsion} we get also the following result.
\begin{lem}\label{Quotient tf}
 Let ${\mathfrak E}$ be a torsion-free Higgs sheaf over $X$ which is strictly Gieseker semistable and let ${\mathfrak F}$ be a Higgs subsheaf of it with $0<{\rm rk\,}{\mathfrak F}<{\rm rk\,}{\mathfrak E}$ and $p_{\mathfrak F}=p_{\mathfrak E}$, then ${\mathfrak Q}={\mathfrak E}/{\mathfrak F}$ is torsion-free.
\end{lem}
\noindent {\it Proof:} Suppose ${\mathfrak Q}$ has torsion ${\mathfrak T}$ and let ${\mathfrak Q}'={\mathfrak Q}/{\mathfrak T}$ and ${\mathfrak F}'$ be the kernel of the morphism ${\mathfrak E}\longrightarrow{\mathfrak Q}'$. Then, from the proof of Proposition \ref{Basic prop. torsion} we see 
that $p_{\mathfrak E}=p_{\mathfrak F}\prec p_{{\mathfrak F}'}$, which is a contradiction, because 
${\mathfrak E}$ is in particular Gieseker semistable.  \;\;Q.E.D.  \\ 

\begin{pro}\label{Key prop. for JH} 
Let ${\mathfrak E}$ be a torsion-free Higgs sheaf over $X$ which is strictly Gieseker semistable and let 
${\mathfrak N}$ be a Higgs subsheaf of it with  $0<{\rm rk\,}{\mathfrak N}<{\rm rk\,}{\mathfrak E}$ and $p_{\mathfrak N}=p_{\mathfrak E}$ and let ${\mathfrak Q}={\mathfrak E}/{\mathfrak N}$. Then, ${\mathfrak N}$ and ${\mathfrak Q}$ are both Gieseker semistable Higgs sheaves and $p_{\mathfrak Q}=p_{\mathfrak E}$.  
 \end{pro} 
 \noindent {\it Proof:}  Let ${\mathfrak E}$ and ${\mathfrak N}$ as in the hypothesis of Proposition \ref{Key prop. for JH}. Since $p_{\mathfrak N}=p_{\mathfrak E}$, the Gieseker semistability of ${\mathfrak N}$ is straightforward.\footnote{Indeed, if there exists a Higgs subsheaf $\tilde{\mathfrak N}$ of ${\mathfrak N}$ with $0<{\rm rk\,}\tilde{\mathfrak N}<{\rm rk\,}{\mathfrak N}$ and  $p_{\mathfrak N}\prec p_{\tilde{\mathfrak N}}$, then $\tilde{\mathfrak N}$ is a Higgs subsheaf of ${\mathfrak E}$ with $p_{\mathfrak E}\prec p_{\tilde{\mathfrak N}}$ and this contradicts the Gieseker semistability of ${\mathfrak E}$.} Now,  from Lemma \ref{Quotient tf} we know that 
 ${\mathfrak Q}$ is torsion-free and we have the following Higgs extension
 \begin{equation}
 \xymatrix{
0 \ar[r]  &  {\mathfrak N} \ar[r]  &   {\mathfrak E} \ar[r]  &   {\mathfrak Q} \ar[r] &  0
}\nonumber
\end{equation} 
with ${\rm rk\,}{\mathfrak Q}\neq 0$. Then, by applying Lemma \ref{Rel. p-rk} to this sequence we obtain also that $p_{\mathfrak Q}=p_{\mathfrak E}$ and therefore the Gieseker semistability of ${\mathfrak Q}$ follows.\footnote{In analogy to the case of Higgs subsheaves, if 
$\tilde{\mathfrak Q}$ is a Higgs quotient of ${\mathfrak Q}$ with $0<{\rm rk\,}\tilde{\mathfrak Q}<{\rm rk\,}{\mathfrak Q}$ and $p_{\tilde{\mathfrak Q}}\prec p_{\mathfrak Q}$, then $\tilde{\mathfrak Q}$ is a Higgs quotient of ${\mathfrak E}$
with $p_{\tilde{\mathfrak Q}}\prec p_{\mathfrak E}$ and this contradicts the Gieseker semistability of ${\mathfrak E}$.}
\;\;Q.E.D.  \\

At this point we establish a theorem involving direct sums of Higgs sheaves, this is an analog of a very 
well known result on stability \cite{Cardona 2} (see \cite{Kobayashi} for a classical version of 
this result on stability).

\begin{thm}\label{property direct sum}
Let ${\mathfrak E}$ and ${\mathfrak E}'$ be two torsion-free Higgs sheaves over $X$. Then ${\mathfrak E}\oplus{\mathfrak E}'$ is Gieseker semistable if and only if ${\mathfrak E}$ and 
${\mathfrak E}'$ are both Gieseker semistable with $p_{\mathfrak E}=p_{\mathfrak{E}'}$.
\end{thm}
\noindent {\it Proof:} Suppose first that ${\mathfrak E}$ and ${\mathfrak E}'$ are both Gieseker semistable with $p=p_{\mathfrak E}=p_{{\mathfrak E}'}$ and let ${\mathfrak F}$ be a Higgs subsheaf of 
${\mathfrak E}\oplus\mathfrak{E}'$ with $0<{\rm rk\,}{\mathfrak F}<{\rm rk}({\mathfrak E}\oplus{\mathfrak E}')$. Then we have the following commutative diagram 
\begin{displaymath}
 \xymatrix{
           &       0                \ar[d]    &     0                    \ar[d]                &     0             \ar[d]         &     \\
0 \ar[r]   &   {\mathfrak S} \ar[r] \ar[d]    &     {\mathfrak F} \ar[r] \ar[d]                &     {\mathfrak Q} \ar[d] \ar[r]  &    0 \\
0 \ar[r]   &   {\mathfrak E} \ar[r]   &     {\mathfrak E}\oplus{\mathfrak E}' \ar[r]    &     {\mathfrak E}' \ar[r]   &    0 \\
}
\end{displaymath}
in which by definition ${\mathfrak S}={\mathfrak F}\cap({\mathfrak E}\oplus 0)$ and ${\mathfrak Q}$ is the image of ${\mathfrak F}$ under 
the projection ${\mathfrak E}\oplus{\mathfrak E}'\longrightarrow{\mathfrak E}'$. In this diagram, the horizontal sequences are 
Higgs extensions of ${\mathfrak F}$ and ${\mathfrak E}\oplus{\mathfrak E}'$. In particular, from the Higgs 
extension of ${\mathfrak E}\oplus{\mathfrak E}'$ and Lemma \ref{Rel. p-rk} it follows that 
$p_{{\mathfrak E}\oplus{\mathfrak E}'}=p$. Now, since ${\mathfrak E}$ and ${\mathfrak E}'$ are both Gieseker 
semistable, we have (for sufficiently large integers $k$)
\begin{equation}
 \chi({\mathfrak S}(k))\le p\,{\rm rk\,}{\mathfrak S}\,, \quad\quad\quad \chi({\mathfrak Q}(k))\le p\,{\rm rk\,}{\mathfrak Q}\,. \nonumber 
\end{equation}
From these inequalities and formula (\ref{Basic prop chi}) applied to the Higgs extension of 
${\mathfrak F}(k)$ (obtained by tensoring the Higgs extension of ${\mathfrak F}$ by ${\mathfrak H}^{k}$) we get 
\begin{equation}
 \chi({\mathfrak F}(k)) =  \chi({\mathfrak S}(k)) + \chi({\mathfrak Q}(k)) \le p\,({\rm rk\,}{\mathfrak S} + {\rm rk\,}{\mathfrak Q}) = p\,{\rm rk\,}{\mathfrak F} \nonumber 
\end{equation}
and hence ${\mathfrak E}\oplus{\mathfrak E}'$ is Gieseker semistable.\\

Conversely, suppose that ${\mathfrak E}\oplus{\mathfrak E}'$ is Gieseker semistable. Clearly, by symmetry we 
can consider also ${\mathfrak E}'$ and ${\mathfrak E}$ as a Higgs subsheaf and a quotient of their direct sum, 
respectively. Then, from the Gieseker semistability of the direct sum and Proposition \ref{Prop. G-Hss} we get 
$p_{{\mathfrak E}\oplus{\mathfrak E}'}=p_{\mathfrak E}=p_{{\mathfrak E}'}$. Now, let ${\mathfrak N}$ be a Higgs subsheaf of ${\mathfrak E}$ with $0<{\rm rk\,}{\mathfrak N}<{\rm rk\,}{\mathfrak E}$, since it is also a Higgs subsheaf of ${\mathfrak E}\oplus{\mathfrak E}'$, then $p_{\mathfrak N}\preceq p_{{\mathfrak E}\oplus{\mathfrak E}'} = p_{\mathfrak E}$ and hence ${\mathfrak E}$ is Gieseker semistable. A similar argument shows the Gieseker semistability of ${\mathfrak E}'$.  \;\;Q.E.D.  \\

At this point we can establish the main result of this paper, this is an extension of a classical result involving morphisms between Gieseker semistable coherent sheaves. As in the classical case, this result on Higgs sheaves also plays an important role in the theory; in fact, some important results can be obtained as a consequence of this theorem. 

\begin{thm}\label{Thm morphisms}
 Let ${\mathfrak E}$ and ${\mathfrak E}'$ be two Gieseker semistable Higgs sheaves over $X$ with 
 ranks $r$ and $r'$, respectively. Let $f:{\mathfrak E}\longrightarrow {\mathfrak E}'$ be a morphism between these 
 Higgs sheaves. Then we have the following: \\
 {\bf (i)} If $p_{{\mathfrak E}'} \prec p_{\mathfrak E}$, then $f=0$ (i.e., it is the zero morphism);\\ 
 {\bf (ii)} If $p_{\mathfrak E} = p_{{\mathfrak E}'}$ and ${\mathfrak E}$ is Gieseker stable, then 
            $r = {\rm rk\,}f({\mathfrak E})$ and $f$ is injective unless $f=0$; \\
  {\bf (iii)} If $p_{\mathfrak E} = p_{{\mathfrak E}'}$ and ${\mathfrak E}'$ is Gieseker stable, then 
            $r' = {\rm rk\,}f({\mathfrak E})$ and $f$ is generically surjective unless $f=0$.          
\end{thm}
\noindent {\it Proof:} Suppose that ${\mathfrak E}$ and ${\mathfrak E}'$ are both Gieseker semistable Higgs sheaves and let ${\mathfrak M}=f({\mathfrak E})$. \\

If ${\rm rk\,}{\mathfrak M}=0$, then necessarily ${\mathfrak M}=0$ 
because ${\mathfrak E}'$ is torsion-free, consequently $f=0$ and there is nothing to prove (all statements follow  trivially). 
Therefore, we assume in the following that ${\rm rk\,}{\mathfrak M}\neq0$. In that case, ${\mathfrak M}$ is a Higgs subsheaf of ${\mathfrak E}'$ and also a Higgs quotient of ${\mathfrak E}$ with $0<{\rm rk\,}{\mathfrak M}\le {\rm min}\{r,r'\}$. Now, if ${\mathfrak M} = {\rm min}\{r,r'\}$ we have two special cases: ${\rm rk\,}{\mathfrak M}=r$ or ${\rm rk\,}{\mathfrak M}=r'$. \\

If ${\rm rk\,}{\mathfrak M}=r$ and ${\mathfrak K}$ denotes the kernel of ${\mathfrak E}\longrightarrow{\mathfrak M}$, 
then we have a Higgs extension
\begin{equation}
 \xymatrix{
0 \ar[r]  &  {\mathfrak K} \ar[r]  &   {\mathfrak E} \ar[r]  &   {\mathfrak M} \ar[r] &  0  
} \nonumber
\end{equation}
with ${\rm rk\,}{\mathfrak K}=0$ and since ${\mathfrak E}$ is torsion-free, necessarily ${\mathfrak K}=0$, so 
${\mathfrak E}$ and ${\mathfrak M}$ are isomorphic and we get $p_{\mathfrak E}=p_{\mathfrak M}$. 
\\

If ${\rm rk\,}{\mathfrak M}=r'$, we have the Higgs extension
\begin{equation}
 \xymatrix{
0 \ar[r]  &  {\mathfrak M} \ar[r]  &   {\mathfrak E}' \ar[r]  &   {\mathfrak T} \ar[r] &  0  
} \nonumber
\end{equation}
in which ${\mathfrak T}$ is torsion. Now, tensoring this sequence by ${\mathfrak H}^{k}$ we obtain 
\begin{equation}
 \xymatrix{
0 \ar[r]  &  {\mathfrak M}(k) \ar[r]  &   {\mathfrak E}'(k) \ar[r]  &   {\mathfrak T}(k) \ar[r] &  0 
} \nonumber
\end{equation}
and, as in the proof of Proposition \ref{Basic prop. torsion}, for a sufficiently large integers $k$ we get 
$\chi({\mathfrak T}(k))>0$. From this fact and (\ref{Basic prop chi}) we obtain that (for sufficiently large integers $k$) \begin{equation}
 \chi({\mathfrak E}'(k)) = \chi({\mathfrak M}(k)) + \chi({\mathfrak T}(k)) > \chi({\mathfrak M}(k)) \nonumber
\end{equation}
holds and hence $p_{\mathfrak M} \prec p_{{\mathfrak E}'}$. \\

To prove (i), let us assume $p_{{\mathfrak E}'}\prec p_{\mathfrak E}$. Since $0<{\rm rk\,}{\mathfrak M}\le {\rm min}\{r,r'\}$, from the above arguments, the hypothesis and Proposition \ref{Prop. G-Hss}, we conclude that for all possible values of  ${\rm rk\,}{\mathfrak M}$ 
\begin{equation}
p_{\mathfrak M} \preceq p_{{\mathfrak E}'} \prec p_{\mathfrak E} \preceq p_{\mathfrak M}\,, \nonumber
\end{equation}
which is clearly impossible at least $f=0$, and hence (i) follows. \\

Let us assume that $p_{\mathfrak E} = p_{{\mathfrak E}'}$ and $f\neq0$ and suppose that ${\mathfrak E}$ is Gieseker stable. If ${\rm rk\,}{\mathfrak M}<r$, then 
\begin{equation}
p_{\mathfrak M} \preceq p_{{\mathfrak E}'} = p_{\mathfrak E} \prec p_{\mathfrak M}\,, \nonumber
\end{equation}
which is impossible and hence $r = {\rm rk\,}{\mathfrak M}$. Moreover, since ${\mathfrak E}$ is torsion-free $f$ is necessarily injective, which proves (ii). \\

If now ${\mathfrak E}'$ is Gieseker stable and ${\rm rk\,}{\mathfrak M}<r'$, then 
\begin{equation}
p_{\mathfrak M} \prec p_{{\mathfrak E}'} = p_{\mathfrak E} \preceq p_{\mathfrak M}\,, \nonumber
\end{equation}
and we have again a contradiction. Therefore $r' = {\rm rk\,}{\mathfrak M}$ and $f$ is 
generically surjective\footnote{Generically surjective means that $f$ is surjective when it is restricted to an open set of 
$X$. In this case, $f$ is clearly surjective when it is restricted to $X\backslash S({\mathfrak M})\cup S({\mathfrak E}')$.} and (iii) follows. \;\;Q.E.D.  \\  

As in the case of coherent sheaves, as a consequence of Theorem \ref{Thm morphisms} we have the following result, its proof is a natural adaptation to the Higgs case of the classical proof of Kobayashi \cite{Kobayashi}. 
\begin{cor}
 Let ${\mathfrak E}$ be a Gieseker stable Higgs sheaf over $X$. Then any endomorphism of ${\mathfrak E}$ is 
 proportional to the identity, i.e., ${\mathfrak E}$ is Higgs simple. 
\end{cor}
\noindent {\it Proof:}  Let $f:{\mathfrak E}\longrightarrow {\mathfrak E}$ be an endomorphism of 
${\mathfrak E}$ and let $a$ be an eigenvalue of $f_{x}:{\mathfrak E}_{x}\longrightarrow {\mathfrak E}_{x}$ at any fixed point $x\in X$. Then $aI$, where $I$ denotes identity endomorphism of the coherent sheaf $E$, is clearly an endomorphism of ${\mathfrak E}$. Then, by applying the part (ii) of Theorem  \ref{Thm morphisms} to the morphism $f - aI$, it follows that  $f - aI$ is injective unless $f - aI = 0$. If 
$f - aI$ is injective, then from Proposition \ref{prop. morph. det}, it induces an injective endomorphism 
${\rm det}(f - aI)$ of the Higgs line bundle ${\rm det\,}{\mathfrak E}$. Now, for this line bundle, such an endomorphism cannot have zeros and consequently $f - aI = 0$. \;\;Q.E.D.  \\  

\begin{cor}
 Let us consider the Higgs extension (\ref{Higgs extension}) over $X$. If ${\mathfrak F}$ and 
 ${\mathfrak Q}$ are both Gieseker semistable with $p_{\mathfrak F}=p_{\mathfrak Q}=p$, then 
 ${\mathfrak E}$ is also Gieseker semistable with $p_{\mathfrak E}=p$. 
\end{cor}
\noindent {\it Proof:} From Lemma \ref{Rel. p-rk} it is clear that $p_{\mathfrak E}=p$. Suppose now that 
${\mathfrak E}$ is not Gieseker semistable, then there exists a Higgs subsheaf ${\mathfrak N}$ of 
${\mathfrak E}$ with $0<{\rm rk\,}{\mathfrak N}<{\rm rk\,}{\mathfrak E}$ and $p \prec p_{\mathfrak N}$. Without loss of generality we can 
assume that ${\mathfrak N}$ is Gieseker semistable.\footnote{Otherwise, we can destabilize 
${\mathfrak N}$ with a Higgs subsheaf ${\mathfrak N}'$. Clearly, this process finishes after a finite number of 
steps, since in the extreme case we get a rank one Higgs subsheaf, which is always Gieseker stable.} Then we get 
the following exact diagram
\begin{displaymath}
 \xymatrix{
            &                           &         0             \ar[d]       &                        &     \\
           &                            &         {\mathfrak N} \ar[d]        &                        &     \\
0 \ar[r]   &   {\mathfrak F} \ar[r]     &          {\mathfrak E} \ar[r]       &  {\mathfrak Q} \ar[r]  &    0 \\
}
\end{displaymath}
and we have a natural morphism $f:{\mathfrak N}\longrightarrow{\mathfrak Q}$ between Gieseker semistable Higgs sheaves 
with $p_{\mathfrak Q} = p\prec p_{\mathfrak N}$. At this point, from part (i) of Theorem \ref{Thm morphisms} 
it follows that $f=0$. Hence, ${\mathfrak N}$ is a Higgs subsheaf of ${\mathfrak F}$ with   
$p_{\mathfrak F} =p \prec p_{\mathfrak N}$, but this contradicts the Gieseker semistability of ${\mathfrak F}$. \;\;Q.E.D.  \\

As we said before, the main result of \cite{Biswas-Schumacher} establishes an equivalence between the notion of polystability and 
the existence of admissible ${\rm HYM}$-metrics for reflexive Higgs sheaves. Now, in the case of locally free Higgs 
sheaves \cite{Bruzzo-Granha, Cardona 1, Jiayu-Zhang}, we also have an equivalence between the notion of semistability and 
the existence of ${\rm apHYM}$-metrics. From these results and Proposition \ref{H-s-->G-H-s} we get the following result.  
\begin{pro}\label{Cor G-H-s vs admissible HYM}
 Let ${\mathfrak E}$ be a reflexive Higgs sheaf over $X$. Then, \\
 {\bf (i)} If ${\mathfrak E}$ has an admissible ${\rm HYM}$-metric, then ${\mathfrak E}=\bigoplus_{i=1}^{s}{\mathfrak E}_{i}$ 
 with each ${\mathfrak E}_{i}$ a Gieseker stable Higgs sheaf;\\
 {\bf (ii)} If moreover, ${\mathfrak E}$ is locally free and it is Gieseker semistable, then there exists 
 an ${\rm apHYM}$-metric on it.
\end{pro}

As it is well known \cite{Bruzzo-Granha}, any locally free Higgs sheaf $\mathfrak E$ over a compact 
K\"ahler manifold $X$ with ${\rm deg\,}{\mathfrak E}<0$ and admitting an apHYM-metric, has no nonzero $\phi$-invariant sections. Then, part (ii) of Proposition \ref{Cor G-H-s vs admissible HYM} immediately implies the following result.
\begin{cor}
 Let ${\mathfrak E}$ be a locally free Higgs sheaf over $X$ with ${\rm deg\,}{\mathfrak E}<0$. If it is Gieseker semistable, then ${\mathfrak E}$ admits no nonzero $\phi$-invariant sections.
\end{cor}

On the other hand, a classical result of Kobayashi \cite{Kobayashi} shows that any holomorphic vector bundle over a compact K\"ahler manifold with an approximate Hermitian-Einstein structure satisfies a Bogomolov-L\"ubke inequality. This result has been extended to Higgs bundles in \cite{Cardona 1}, and hence we have a Bogomolov-L\"ubke inequality for any Higgs bundle admitting an apHYM-metric; as a consequence of this fact and part (ii) of Proposition \ref{Cor G-H-s vs admissible HYM}, we get the following result for locally free Higgs sheaves.
\begin{cor}\label{Bogomolov ineq.}
Let ${\mathfrak E}$ be a locally free Higgs sheaf over $X$. If ${\mathfrak E}$ is Gieseker semistable, then
\begin{equation}
\int_{X}[2r c_{2}(\mathfrak E) - (r-1)c_{1}(\mathfrak E)^{2}]\cdot\omega_{H}^{n-2} \ge 0\,. \nonumber
\end{equation}
\end{cor}

\section{Final remarks}

Let ${\mathfrak E}$ be a Gieseker semistable Higgs sheaf over $X$. Following the classical definition for coherent sheaves \cite{Gieseker, Huybrechts-Lehn, Kobayashi} and in analogy to the definition of flags for Higgs sheaves \cite{Cardona 4}, a Jordan-H\"older filtration of ${\mathfrak E}$ is a family 
$\{{\mathfrak E}_{i}\}_{i=0}^{s+1}$ of Higgs subsheaves of ${\mathfrak E}$ with
\begin{equation}
0 = {\mathfrak E}_{s+1}\subset{\mathfrak E}_{s}\subset\cdots\subset{\mathfrak E}_{1}\subset{\mathfrak E}_{0}={\mathfrak E}
\end{equation}
and such that the Higgs quotients ${\mathfrak Q}_{i}={\mathfrak E}_{i}/{\mathfrak E}_{i+1}$ are Gieseker stable and 
$p_{{\mathfrak Q}_{i}} = p_{{\mathfrak E}}$ for $i=0,...,s$. With the quotients ${\mathfrak Q}_{i}$ we can define a Higgs sheaf 
\begin{equation}
\mathfrak{gr}({\mathfrak E}) = \bigoplus_{i=0}^{s}{\mathfrak Q}_{i}
\end{equation}
which is commonly called the associated grading to ${\mathfrak E}$. Again, in analogy to the classical literature \cite{Huybrechts-Lehn}, we say that two Gieseker semistable Higgs sheaves ${\mathfrak E}_{1}$ and ${\mathfrak E}_{2}$ over $X$ with $p_{{\mathfrak E}_{1}}=p_{{\mathfrak E}_{2}}$ are  Jordan-H\"older equivalent or S-equivalent, if $\mathfrak{gr}({\mathfrak E}_{1}) \cong \mathfrak{gr}({\mathfrak E}_{2})$. The construction of a Jordan-H\"older filtration for Gieseker semistable Higgs sheaves is similar to the classical case. Suppose that ${\mathfrak E}$ is a Gieseker semistable Higgs sheaf and let ${\mathfrak E}_{0}={\mathfrak E}$. If  ${\mathfrak E}_{0}$ is also Gieseker stable, then we take ${\mathfrak E}_{1}=0$ and we have a Jordan-H\"older filtration with $s=0$. If it is not, i.e., if ${\mathfrak E}_{0}$ is strictly Gieseker semistable, we consider the set of all Higgs subsheaves of ${\mathfrak E}$ satisfying the hypothesis of Proposition \ref{Key prop. for JH}, then we take from this set a Higgs subsheaf of maximal rank and we denote this as ${\mathfrak E}_{1}$. \\
 
From Proposition \ref{Key prop. for JH} we know that the Higgs quotient ${\mathfrak Q}_{0}={\mathfrak E}_{0}/{\mathfrak E}_{1}$ is Gieseker semistable with $p_{{\mathfrak Q}_{0}}=p_{\mathfrak E}$. Clearly, among all Higgs quotients of ${\mathfrak E}$ satisfying these conditions, ${\mathfrak Q}_{0}$ is one with minimal rank and therefore, it is also Gieseker stable. In fact, if there exists a Higgs quotient $\tilde{\mathfrak Q}$ of ${\mathfrak Q}_{0}$ with $0<{\rm rk\,}\tilde{\mathfrak Q}<{\rm rk\,}{\mathfrak Q}_{0}$ and $p_{\tilde{\mathfrak Q}}\preceq p_{{\mathfrak Q}_{0}}$, since 
$\tilde{\mathfrak Q}$ is also a Higgs quotient of $\mathfrak E$ and it is Gieseker semistable, then $p_{\mathfrak E}\preceq p_{\tilde{\mathfrak Q}}$ and hence, necessarily $p_{\tilde{\mathfrak Q}} = p_{\mathfrak E}$. But this contradicts the fact that ${\mathfrak Q}_{0}$ has minimal rank. Now, if ${\mathfrak E}_{1}$ is Gieseker stable, we take ${\mathfrak E}_{2}=0$ and we have a Jordan-H\"older filtration of ${\mathfrak E}$ with $s=1$. Otherwise, we consider the same procedure with ${\mathfrak E}_{1}$ instead of 
 ${\mathfrak E}_{0}$, in this way we get always Gieseker stable Higgs quotients. Clearly, after a finite number of steps this procedure finishes with a stable Higgs sheaf (in the extreme case we obtain a Higgs sheaf or rank one which is Gieseker stable) and this proves the existence of a Jordan-H\"older filtration. \\
 
From the construction of the Jordan-H\"older filtration and Proposition \ref{Key prop. for JH} it follows that the Higgs sheaves 
${\mathfrak E}_{i}$ with $i=0,...,s$ of the Jordan-H\"older filtration are Gieseker semistable with $p_{{\mathfrak E}_{i}}=p_{\mathfrak E}$. Moreover, ${\mathfrak E}_{s}$ is always Gieseker stable. This in part shows that a Higgs sheaf may have different Jordan-H\"older filtrations. In fact, if ${\mathfrak L}$ and ${\mathfrak L}'$ are two Higgs line bundles over $X$ with $p_{\mathfrak L}=p_{{\mathfrak L}'}$, they are Gisesker stable and by Theorem \ref{property direct sum} we know that the direct sum ${\mathfrak L}\oplus{\mathfrak L}'$ is Gieseker semistable. Clearly, ${\mathfrak L}$ and ${\mathfrak L}'$ define two different Jordan-H\"older filtrations of ${\mathfrak L}\oplus{\mathfrak L}'$. Although the Jordan-H\"older filtration of a Higgs sheaf is in general not unique, its associated grading is. Indeed (see Simpson \cite{Simpson M1, Simpson M2} for more details) we have the following result. 
\begin{thm}\label{Existence JH}
Let ${\mathfrak E}$ be a torsion-free Higgs sheaf over $X$. If it is Gieseker semistable, then 
${\mathfrak E}$ has a Jordan-H\"older filtration. Furthermore, up to isomorphism, the Higgs sheaf $\mathfrak{gr}({\mathfrak E})$ does not depend on the choice of the Jordan-H\"older filtration.
\end{thm}
 
Finally, it is important to mention that there exists another important kind of filtrations for coherent sheaves known as Harder-Narasimhan filtrations \cite{Huybrechts-Lehn, Kobayashi}; a similar definition can be done in the Higgs case. Let ${\mathfrak E}$ be a Higgs sheaf over $X$, a Harder-Narasimhan filtration of ${\mathfrak E}$ is a family 
$\{{\mathfrak E}_{i}\}_{i=0}^{l}$ of Higgs subsheaves of $\mathfrak E$ with
\begin{equation}
0 = {\mathfrak E}_{0}\subset{\mathfrak E}_{1}\subset\cdots\subset{\mathfrak E}_{l-1}\subset{\mathfrak E}_{l}={\mathfrak E}
\end{equation}
and such that the Higgs quotients ${\mathfrak Q}_{i}={\mathfrak E}_{i}/{\mathfrak E}_{i-1}$ are all Gieseker semistable with
$p_{{\mathfrak Q}_{i}} >p_{{\mathfrak Q}_{i+1}}$ for $i=1,...,l-1$. Harder-Narasimhan filtrations do exist for Higgs sheaves and (in contrast to Jordan-H\"older filtrations) they are unique. The existence of Harder-Narasimhan filtrations for Higgs sheaves can be seen as a particular case of a more general approach \cite{Simpson M1, Simpson M2}. In the context of Higgs sheaves this result can be written as follows.
\begin{thm}\label{Existence HN}
Let ${\mathfrak E}$ be a torsion-free Higgs sheaf over $X$. Then, ${\mathfrak E}$ has a unique Harder-Narasimhan filtration.
\end{thm}

\noindent{\bf Aknowledgements}\\

\noindent Part of this paper was done when both authors were postdoctoral fellows at Centro de Investigaci\'on en Matem\'aticas (CIMAT) in Guanajuato, M\'exico. The authors would like to thank CIMAT and CONACyT for support. The first author wants to thank also U. Bruzzo for some important comments and for drawing Simpson's papers \cite{Simpson M1, Simpson M2} to their attention. Finally, the second author wants to thank the PRO-SNI programe of the Universidad de Guadalajara for partial support.

{\sc Conacyt Research Fellow - Instituto de Matem\'aticas, Universidad Nacional Aut\'onoma de 
M\'exico - UNAM, Leon 2 altos, Col. centro,  C.P. 68000, Oaxaca de Ju\'arez, Oax., M\'exico}.\\
{\it E-mail addresses}: {\bf sholguin@im.unam.mx , sholguinca@conacyt.mx} \\

{\sc Departamento de Matem\'aticas - CUCEI, Universidad de Guadalajara., Av. Revoluci\'on 1500,   
C.P. 44430, Guadalajara Jal., M\'exico}.\\
{\it E-mail addresses}: {\bf osbaldo@cimat.mx , osbaldo.mata@academico.udg.mx}

\end{document}